\documentclass[11pt]{article}

\usepackage[T1]{fontenc}
\usepackage[utf8]{inputenc}
\usepackage{lmodern}
\usepackage{amsmath,amssymb,amsfonts,amsthm}
\allowdisplaybreaks
\usepackage{mathtools}
\usepackage{mathrsfs}
\usepackage{bm}
\usepackage{graphicx}
\usepackage{tikz}
\usepackage{pgfplots}
\pgfplotsset{compat=1.16}
\usepackage{booktabs,tabularx}
\usepackage{enumitem}

\usepackage[colorlinks=true,
linkcolor=blue,
citecolor=blue,
urlcolor=blue]{hyperref}

\usepackage[a4paper,margin=3cm]{geometry}

\newtheorem{theorem}{Theorem}[section]
\newtheorem{proposition}[theorem]{Proposition}
\newtheorem{lemma}[theorem]{Lemma}
\newtheorem{corollary}[theorem]{Corollary}

\theoremstyle{definition}

\newtheorem{remark}[theorem]{Remark}

\newcommand{\R}{\mathbb{R}}

\newcommand{\bh}{b_h}
\newcommand{\q}{q^{\mathrm{Rob}}_{m,h}}
\newcommand{\NRob}{N}
\newcommand{\dd}{\,\mathrm{d}}

\title{Counting eigenvalues below the lowest Landau level: Robin–Dirichlet boundary conditions in an annulus}

\author{
Marwa Nasrallah$^{1,2,\ast}$, Iman Al Moustafa$^1$
}

\date{}

\begin{document}
\maketitle 

\begin{center}
{\small
$^{1}$ Department of Mathematics and Physics, School of Arts and Sciences,
Lebanese International University, Bekaa, Lebanon\\

$^{2}$ Department of Mathematics, Faculty of Sciences IV, 
Lebanese University, Zahle, Lebanon\\

$^\ast$  \texttt{marwa.nasrallah@liu.edu.lb}}

\end{center}

\begin{abstract}
We study the magnetic Laplacian with constant field in a planar annulus, in the semiclassical limit, with a Dirichlet condition on the inner boundary and a Robin condition on the outer boundary whose parameter $\beta(h)$ converges to a finite limit $\beta$. We determine the leading asymptotics of the number of eigenvalues below the lowest Landau level for every $\beta\in\mathbb{R}$. Our main result concerns the attractive regime $\beta<0$. There, the outer boundary binds states at energy of order one below the lowest Landau level, and these states compete with those localized in the interior. The resulting limiting density is continuous and piecewise linear in $\beta$, with a kink at an explicit critical value $-\beta_c$ where its slope jumps from $-2$ to $-4$. For $-\beta_c\le\beta<0$ the attractive boundary increases the count; for $\beta<-\beta_c$ the count is carried by the boundary states and no longer depends on the inner radius. Since these states are invisible to the harmonic approximation, we control them by a sharp semiclassical trace inequality with optimal constant, combined with an Agmon estimate whose weight is frozen in the boundary layer. In the repulsive regime $\beta>0$ the count decreases linearly and becomes negligible beyond $\beta_c$.
\end{abstract}


\section{Introduction and main result}

The number of eigenvalues of a magnetic Laplacian below the lowest Landau level measures how strongly a boundary can bind charged particles in a magnetic field. For a constant field in two dimensions, the lowest Landau level is located at energy $h$. In a bounded domain the spectrum below this threshold comes from states localized near the boundary, and the boundary condition decides whether such states exist. Under a Dirichlet condition the spectrum lies above $h$, whereas under a Neumann condition a family of edge states propagating along the boundary appears below it \cite{HelfferMorame, FournaisHelffer, RaymondEdge}; see \cite{RaymondAnnulus} for a general account and \cite{Landau} for the physical background.

In this paper we show that an \emph{attractive} Robin condition changes this picture qualitatively. It creates a second, genuinely different family of boundary states, bound at energy of order one below the lowest Landau level, and the competition between this family and the states localized in the interior produces a sharp transition in the leading-order count at an explicit critical value of the Robin parameter. We establish this in an annulus, with a Dirichlet condition on the inner boundary, where the rotational symmetry allows a complete description for every value of the Robin parameter.

Fix $R \in (0,1)$ and define the annulus 
\begin{equation}\label{Annulus}
\Omega_R = \{ x \in \mathbb{R}^2 : R < |x| < 1 \}.
\end{equation}
Consider the magnetic Laplacian
\begin{equation}\label{eq:Ope}
P_h = -(h \nabla -i{ A})^2 \quad \text{on } \Omega_R, 
\end{equation}
where
\begin{equation}\label{MagPot}
{A}(x)=\frac12(-x_2,x_1).
\end{equation}
The magnetic field is constant, since
\[
\operatorname{curl} A = 1.
\]
Throughout this work, for $A = (A_1, A_2) \in C^1(\Omega; \R^2)$, we define
\(
\operatorname{curl} A:= \partial_{x_1} A_2-\partial_{x_2} A_1.
\)
The parameter $h>0$ denotes the semi-classical parameter and will be taken small, $h\rightarrow 0$. 
Impose the following boundary conditions
\begin{itemize}
\item Dirichlet condition on the inner boundary $|x| = R: u = 0,$
\item Robin condition on the outer boundary $|x| = 1:
\nu \cdot (h\nabla - iA)u + \beta(h) u = 0.$
\end{itemize}
Here $\nu$ is the outward unit normal to the boundary.
We assume that
\[
\beta(h)\to\beta,
\qquad h\to0.
\]

We study the semiclassical limit $h \to 0$ and define
\begin{equation}\label{NB-OF-EIG}
N(h)
=
\#\{\lambda \in \mathrm{Spec}(P_h) : \lambda < h\},
\end{equation}
that is, the number of eigenvalues of $P_h$ below the lowest Landau level.

The regime $\beta(h)\to\beta$ with a finite limit is the critical one for the
counting problem. By the exact identity of Proposition~\ref{ExactIdentity}, the
boundary coefficient $2mh-1+2\beta(h)$ changes sign at $2mh=1-2\beta(h)$, so a
finite limit $\beta\ne0$ moves the range of contributing angular momenta by an
amount of order $h^{-1}$, which is the order of the count itself. If
$\beta(h)\to0$ the leading term is that of the Neumann problem; if
$\beta(h)\to+\infty$ the count is $o(h^{-1})$, and if $\beta(h)\to-\infty$ one expects the boundary states to dominate (this case is not treated here). Physically, $\beta$ models a surface admittance or impedance
at the outer boundary, as in the de~Gennes boundary condition of surface
superconductivity (see \cite{CorreggiRougerie, FournaisHelfferBook}).

\subsection{Main result} 
The main result gives the leading-order asymptotics of the counting function $N(h)$ for every value of the limiting Robin parameter $\beta$.
\bigskip\begin{theorem}[\bf Eigenvalue counting with Robin boundary condition]\label{thm:main_Robin}
Assume that $\beta(h)\to\beta\in\mathbb{R}$ as $h\to0$, and let
\begin{equation}\label{Def-Re}
R_e=\sqrt{\dfrac{1-R^2}{2|\ln R|}}.
\end{equation} Then for every $\beta\in\mathbb R$,
\begin{equation}\label{eq:ell-all}
\lim_{h\to0}2h\,N(h)=\ell(\beta):=\max\Bigl\{\bigl(1-R_e^2-2\beta\bigr)_+,\;-4\beta\Bigr\}.
\end{equation}
\end{theorem}

\subsection{The attractive regime and the transition at \texorpdfstring{$-\beta_c$}{-beta\_c}}

Set $\beta_c:=\frac{1-R_e^2}{2}$. By Theorem~\ref{thm:main_Robin}, the limiting density $\ell(\beta)$ is continuous and piecewise linear, with slope $-4$ on $(-\infty,-\beta_c)$, slope $-2$ on $(-\beta_c,\beta_c)$ and $\ell\equiv0$ on $[\beta_c,\infty)$; see Figure~\ref{fig:density-vs-beta}. The new phenomena occur for $\beta<0$.
\begin{itemize}
\item \emph{Boundary states.} An attractive Robin condition creates, for every angular momentum $m$ with $|1-2mh|<2|\beta|$, a fiber eigenvalue at distance of order one below the lowest Landau level (Proposition~\ref{prop:bs}). These states live in a layer of width $\mathcal{O}(h)$ at the outer boundary, and they have no analogue for Dirichlet, Neumann or repulsive Robin conditions.
\item \emph{Enhancement of the count.} For $-\beta_c\le\beta<0$ the density $1-R_e^2-2\beta$ exceeds the Neumann value $1-R_e^2$: an attractive boundary increases the number of eigenvalues below $h$ at leading order.
\item \emph{Transition at $-\beta_c$.} At $\beta=-\beta_c$ the density has a kink, where its slope jumps from $-2$ to $-4$. For $\beta<-\beta_c$ the count is carried by the boundary states, the density equals $-4\beta$, and it no longer depends on the inner radius $R$. The critical value is explicit and is fixed by the effective radius $R_e$, that is, by the tunnelling competition between the inner Dirichlet layer and the outer boundary.
\end{itemize}

\subsection{Main difficulties}

The boundary states are invisible to the harmonic approximation near the minimum $r_*$ of the effective potential, which is the standard tool for counting eigenvalues below the Landau level. They cannot be controlled by the usual trace inequality
$|u(1)|^2\le\varepsilon\int_R^1|u'|^2r\,dr+C_\varepsilon\int_R^1|u|^2r\,dr$ either: with $\varepsilon\sim h$ it loses a term of size $\mathcal{O}(1)\|u\|^2$, while the relevant gap is of order $h$, and this loss is exactly the energy $-\beta^2$ of a boundary state. We introduce two tools to overcome this. The first is a semiclassical trace inequality whose constant $b^2$ is optimal (Lemma~\ref{lem:trace}), so that the boundary-state energy is captured without loss. The second is an Agmon estimate with a weight frozen in the boundary layer (Proposition~\ref{prop:agmon}), which separates the boundary states from the states localized in the interior well. Together they give the leading-order count for every $\beta<0$ (Section~\ref{sec:neg}).

For $\beta\ge0$, the Robin term shifts the upper threshold of contributing angular momenta from $2mh<1$ to $2mh<1-2\beta$. Unlike the Neumann term, it does not vanish on the boundary trace of test functions, so the quasimode construction for favorable angular momenta (Proposition~\ref{FavAngMom}) requires an explicit control of the boundary contribution $\beta(h)h|u(1)|^2$ against the exponentially small tail of the quasimode; see \eqref{QF-h}--\eqref{equation-inv-rstar}.

\subsection{Relation to previous work}

For $\beta=0$, Theorem~\ref{thm:main_Robin} gives the Dirichlet--Neumann count of Fournais and Kachmar~\cite{FournaisKachmar} (Remark~\ref{rem:neumann}), whose fiber decomposition and Agmon analysis of the inner Dirichlet layer we also use (Section~\ref{ssec:unfav}). Robin and de~Gennes conditions for magnetic Laplacians have been studied for the ground state energy~\cite{Ka06} and for Weyl-type asymptotics~\cite{Ka08}; for the Robin Laplacian without magnetic field see \cite{KKR16, HKR17, PankrashkinPopoff}.

Fournais, Frank, Goffeng, Kachmar and Sundqvist \cite{FFGK} recently
counted the negative eigenvalues of the magnetic Pauli operator in bounded,
possibly multiply connected, planar domains with Robin boundary conditions. Their bound is of Aharonov--Casher
type and is obtained from a decomposition of the quadratic form into a
$\partial_{\bar z}$-part and a boundary term, combined with an index
argument. In the semiclassical limit it yields, for a bounded domain $\Omega$,
$N(h)\sim|\Omega|/(2\pi h)$ eigenvalues of $-(h\nabla-i A)^2$ below
$h$ when the Neumann condition holds on the whole boundary. For the annulus
this corresponds to the density $1-R^2$: every angular momentum with
$r_*\in(R,1)$ contributes. The present count arises from a different,
boundary-driven mechanism. Our quasimodes $r^me^{-r^2/(4h)}$ are of the same
holomorphic type, but the Dirichlet condition at $r=R$ removes the angular
momenta whose well lies below the effective radius $R_e$ through tunnelling
competition between the two boundary layers, and the Robin parameter moves the upper threshold from $2mh<1$ to $2mh<1-2\beta$. Consequently the density is
$(1-R_e^2-2\beta)_+$ instead of $1-R_e^2$ whenever $\beta\ge-\beta_c$.
For $\beta<-\beta_c$ the count is dominated by the boundary states created by
the attractive Robin condition; this contribution, of density $-4\beta$,
has no counterpart in~\cite{FFGK}.

The paper is organized as follows. Section~\ref{FiberDec} uses the rotational symmetry to reduce the problem to a family of one-dimensional operators, and Section~\ref{SectionReductionFirstBand} shows that each of them contributes at most one eigenvalue below the lowest Landau level. Section~\ref{DecFavUnfav} analyses the favorable and unfavorable angular momenta, and Section~\ref{sec:proofmain} completes the proof for $\beta\ge0$. Section~\ref{sec:neg}, the core of the paper, treats the attractive case $\beta<0$: the boundary states, the sharp trace inequality and the Agmon estimate in the boundary layer. Section~\ref{sec:perspectives} discusses open problems.


\subsection{Remarks on the theorem}

\begin{remark}[\bf Strategy of the proof]\label{rem:comparison}
The proof combines the fiber decomposition associated with the rotational symmetry, quasimodes for favorable angular momenta, and Agmon estimates for unfavorable ones. The Robin condition enters in two different ways, depending on the sign of $\beta$.

For $\beta \ge 0$, it enters through the boundary coefficient $2mh - 1 + 2\beta(h)$ in the exact identity~\eqref{ExactIdentityEq}. This coefficient moves the threshold of admissible angular momenta from $2mh < 1$ to $2mh < 1 - 2\beta$, which gives the favorable interval~\eqref{FavInt} and the quasimode estimate~\eqref{QF-h}--\eqref{equation-inv-rstar}. The remaining steps rely on the comparison $q^{\mathrm{Rob}}_{m,h} \ge q^{\mathrm{DN}}_{m,h}$ in~\eqref{eq:compare}.

For $\beta < 0$, this comparison fails. The attractive boundary creates states at energy $\mathcal{O}(1)$ below $h$ (Proposition~\ref{prop:bs}), which the harmonic approximation does not detect. They are controlled by the sharp trace inequality of Lemma~\ref{lem:trace}, whose constant $b^2$ is optimal, and by an Agmon weight frozen in the boundary layer (Proposition~\ref{prop:agmon}). This is the main new part of the paper. For $\beta < -\beta_c$ these states change the leading-order count.
\end{remark}

\begin{remark}[{\bf Interpretation of the effective radius}]
The quantity $R_e$ defined in \eqref{Def-Re} arises from the competition between the inner Dirichlet boundary and the outer boundary layer. It corresponds to the critical value of the minimum point $r_\ast$ of the effective potential  defined below in \eqref{rstar} for which the two boundary contributions are of the same exponential order.
\end{remark}

\begin{remark}[{\bf Reduction to the Neumann case}]\label{rem:neumann}
If $\beta = 0$, the Robin condition reduces to a Neumann condition at leading order, and Theorem~\ref{thm:main_Robin} gives
\[
N(h) \sim \frac{1 - R_e^2}{2h},
\]
which coincides with the Dirichlet--Neumann annulus result in \cite{FournaisKachmar}.
\end{remark}

\begin{figure}[ht]
\centering
\begin{tikzpicture}[x=1.6cm,y=0.75cm,>=stealth]
  \draw[->] (-3.2,0) -- (2.2,0) node[right] {$\beta$};
  \draw[->] (0,0) -- (0,7.4) node[above] {$\ell(\beta)$};
  \draw[very thick,blue!70!black] (-1.75,7) -- (-1,4) -- (1,0) -- (2.1,0);
  \draw[thick,dashed,red!70!black] (-1,4) -- (-1.75,5.5);
  \node[red!70!black,left,align=right,font=\small] at (-1.75,5.5)
     {$1-R_e^2-2\beta$};
  \node[blue!70!black,right,font=\small] at (-1.65,7) {$-4\beta$};
  \draw[densely dotted] (-1,0) node[below] {$-\beta_c$} -- (-1,4);
  \draw[densely dotted] (-1,4) -- (0,4) node[right,font=\small] {$2(1-R_e^2)$};
  \node[below] at (1,0) {$\beta_c$};
  \node[right,font=\small] at (0,2) {$1-R_e^2$};
  \fill (0,2) circle (1.5pt);
\end{tikzpicture}
\caption{The limiting density $\ell(\beta)=\lim_{h\to0}2h\,N(h)$ of
Theorem~\ref{thm:main_Robin}, see \eqref{eq:ell-all}. It equals the Neumann value
$1-R_e^2$ at $\beta=0$ (Remark~\ref{rem:neumann}) and vanishes for
$\beta\ge\beta_c:=(1-R_e^2)/2$ (Remark~\ref{rem:vanishing}). 
 It has slope $-2$ on$(-\beta_c,\beta_c)$ and slope $-4$ on $(-\infty,-\beta_c)$, where boundary states
make it strictly larger than the line $1-R_e^2-2\beta$ (dashed), which gives the density only for $\beta\ge -\beta_c$
(Remark~\ref{rem:negative}).}
\label{fig:density-vs-beta}
\end{figure}
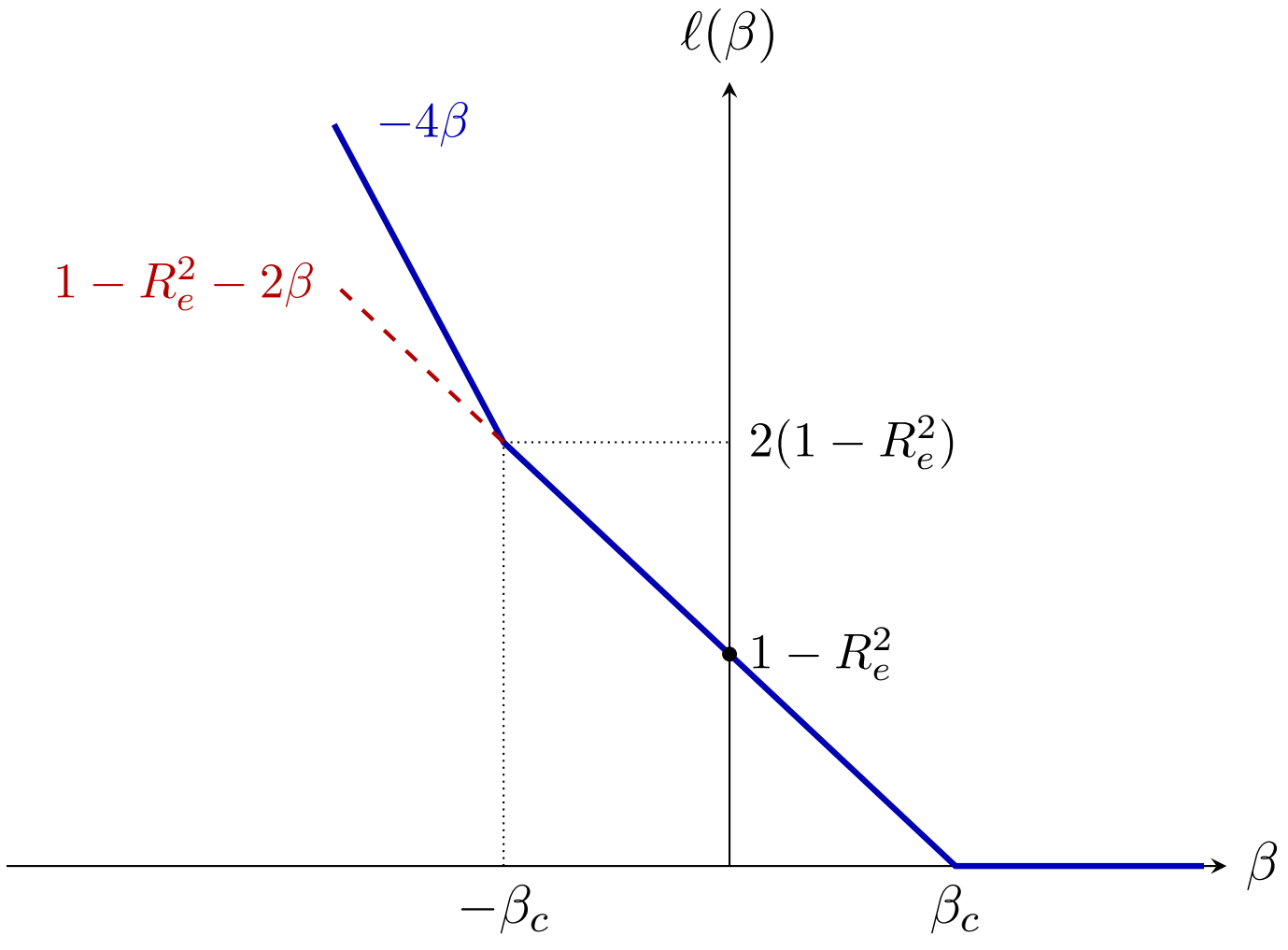

\begin{remark}[{\bf Critical value of the Robin parameter}]\label{rem:vanishing}
If $\beta\ge\beta_c$, then Theorem~\ref{thm:main_Robin} gives
\[
  N(h)=o(h^{-1}),\qquad h\to 0 .
\]
Thus the number of eigenvalues below the lowest Landau level is negligible
compared with the order $h^{-1}$ of the Neumann case $\beta=0$
(Remark~\ref{rem:neumann}). This does not claim that the spectrum below $h$ is empty:
the theorem gives no information on the lower-order term, and finitely many
(or $o(h^{-1})$ many) eigenvalues may remain. For $0\le\beta<\beta_c$ the
leading-order density $1-R_e^2-2\beta$ decreases linearly in $\beta$ and
vanishes exactly at $\beta=\beta_c$ (see Figure~\ref{fig:density-vs-beta}).
\end{remark}

\begin{remark}[{\bf The case $\beta<0$}]\label{rem:negative}
This case is treated in Section~\ref{sec:neg}. \!The
proof for~~$\!\beta\ge0$ does not extend directly: the comparison
$q^{\mathrm{Rob}}_{m,h}\ge q^{\mathrm{DN}}_{m,h}$ fails, and the standard trace inequality
$|u(1)|^2\le\varepsilon\int_R^1|u'|^2r\,dr+C_\varepsilon\int_R^1|u|^2r\,dr$ controls the
attractive boundary term only up to an $O(1)\|u\|^2$ loss (one needs $\varepsilon\sim h$), while the relevant gap is of order $h$. This loss cannot be
avoided: it is the energy $-\beta^2$ of a boundary state. It has to be compared
with the potential at the wall, $V_m(1)=(1-2mh)^2/4$, which is also of order
one. For $-\beta_c\le\beta<0$ the density is still $1-R_e^2-2\beta$, which now exceeds the
Neumann value $1-R_e^2$. For $\beta<-\beta_c$ the boundary states contribute at leading
order and the density becomes $-4\beta$, which is strictly larger than $1-R_e^2-2\beta$.
This is why $\ell(\beta)$ in \eqref{eq:ell-all} is the maximum of the two expressions
(see Figure~\ref{fig:density-vs-beta}).
\end{remark}

\section{Fiber decomposition}\label{FiberDec}





The operator $P_h$ introduced in \eqref{eq:Ope} is the self-adjoint operator associated with the quadratic form
\begin{equation}\label{eq:form}
  q_h(u)=\int_{\Omega_R}\bigl|(h\nabla-iA)u\bigr|^2\dd x
  +\beta(h)\,h\int_{|x|=1}|u|^2\dd s,
\end{equation}
with form domain $\{u\in H^1(\Omega_R):\ u=0 \text{ on } |x|=R\}$.
Its operator domain is determined by the boundary condition
\begin{equation}\label{eq:robin}
  \nu\cdot(h\nabla-iA)u+\beta(h)\,u=0
  \quad\text{on } |x|=1 .
\end{equation}

Since the coefficients of \(P_h\) are independent of \(\theta\), we expand functions into Fourier series,

\[
u(r,\theta)=
\sum_{m\in\mathbb Z}
u_m(r)e^{im\theta}.
\]

This yields the orthogonal decomposition

\[
L^2(\Omega_R)=
\bigoplus_{m\in\mathbb Z}
L^2((R,1),rdr).
\]

Consequently, the operator \(P_h\) admits the orthogonal decomposition
\begin{equation}\label{UE}
P_h
\simeq
\bigoplus_{m\in\mathbb Z}
H_{m,h},
\end{equation}

where

\begin{equation}\label{Hmh}
H_{m,h}=
-h^2\frac{d^2}{dr^2}
-\frac{h^2}{r}\frac{d}{dr}
+
\left(
\frac{mh}{r}-\frac r2
\right)^2 .
\end{equation}
Since $A$ is tangential, $\nu\cdot(h\nabla-iA)u=h\,\partial_r u$ on $|x|=1$.
Hence, on the $m$-th Fourier mode, \eqref{eq:robin} reads
\begin{equation}\label{eq:fiberrobin}
  h\,u_m'(1)+\beta(h)\,u_m(1)=0,
\end{equation}
and the fiber form is
\begin{equation}\label{eq:fiberform}
  q^{\mathrm{Rob}}_{m,h}(u)=\int_R^1\bigl(h^2|u'|^2+V_m|u|^2\bigr)r dr
  +\beta(h)\,h\,|u(1)|^2,
  \qquad V_m(r)=\Bigl(\frac{mh}{r}-\frac r2\Bigr)^2 ,
\end{equation}
with domain $\{u\in H^1(R,1):u(R)=0\}$.

The domain of \(H_{m,h}\) is then determined by the Dirichlet condition on the inner boundary
\(
u(R)=0
\)
and the Robin condition on the outer boundary \(hu^\prime(1)+\beta(h)u(1)=0\).

Since the operator $H_{m,h}$ is an operator with compact resolvent, its spectrum is purely discrete and we denote by
\[
\lambda_0(H_{m,h})
<
\lambda_1(H_{m,h})
\leq
\lambda_2(H_{m,h})
\leq \cdots
\]
its increasing sequence of eigenvalues. The first eigenvalue of $H_{m,h}$ is simple.
It follows from \eqref{UE} that
\[
{\rm Spec} (P_{h})=\underset{m\in\mathbb{Z}}{\bigcup} {\rm Spec} (H_{m,h})=\underset{m\in\mathbb{Z}}{\bigcup} \{\lambda_0(H_{m,h}),\lambda_1(H_{m,h}),\cdots\}
\]
So the quantity $N(h)$ defined in \eqref{NB-OF-EIG} becomes
\begin{equation}\label{NB-OF-EIG1}
N(h)=\#\left\{ (m,j) \in \mathbb{Z}\times\mathbb{N} : \lambda_j(H_{m,h}) < h \right\}.
\end{equation}

\section{Higher order eigenvalues}\label{SectionReductionFirstBand}
For $m\in\mathbb{Z}$, we define
\begin{equation}\label{u}
u_m(r)=r^m e^{-\frac{r^2}{4h}}.
\end{equation}
Since $h\,u_m'/u_m=\frac{mh}{r}-\frac r2$, we have $h^2(u_m'/u_m)^2=V_m$, and a direct computation from \eqref{Hmh} gives
\begin{equation}\label{um-eq}
H_{m,h}u_m=h\,u_m\quad\text{on }(R,1).
\end{equation}
We introduce
\begin{equation}\label{rstar}
r_* = \sqrt{2mh},
\end{equation}
which corresponds to the minimum of the effective potential.
Rewrite $u_m(r)$ as
\begin{equation}\label{Def2-u}
u_m(r)=e^{-\frac{\phi(r)}{4h}},
\end{equation}
where \begin{equation}\label{Defphi}
\phi(r)=r^2-2r_\ast^2 \ln r.
\end{equation}
By computing the first and second derivatives of $\phi$, we obtain
 \begin{equation}\label{FS-Der}
 \phi^\prime(r)=\frac{2r^2-2{r_\ast^2}}{r} \qquad {\rm and}\qquad \phi^{\prime\prime}(r)=2+2\frac{r_\ast^2}{r^2}.
 \end{equation}
 In particular, we find that
\[
\min \phi(r)= \phi(r_\ast)\qquad {\rm and}\qquad \phi^{\prime\prime}(r_\ast)=4.
\]
By Laplace's method, since $\phi''(r_\ast)=4$, for $r_\ast\in[R+\varepsilon,1-\varepsilon]$ we have
\begin{equation}\label{LaplaceMethod}
\begin{aligned}
\int_R^1 |u_m(r)|^2\, r\,dr&= \int_R^1 e^{ -\frac{\phi(r)}{2h}}\, r\,dr
=\bigl(1+o(1)\bigr)\, r_\ast\, e^{ -\frac{\phi(r_\ast)}{2h}}\sqrt{\frac{2\pi\cdot 2h}{\phi''(r_\ast)}}\\
&=\bigl(1+o(1)\bigr)\sqrt{\pi h}\; r_\ast\, e^{ -\frac{\phi(r_\ast)}{2h}},\qquad h\to0.
\end{aligned}
\end{equation}
The fiber decomposition reduces the spectral problem to the study of the family of
one-dimensional operators $H_{m,h}$.
The behavior of $\lambda_0(H_{m,h})$ is governed by the effective potential
\[
V_m(r) = \left(\frac{mh}{r} - \frac{r}{2}\right)^2,
\]
whose minimum is attained at $r_\ast$.
Thus, the position of $r_\ast$ relative to the interval $(R,1)$ determines the nature of
the corresponding eigenfunction.

As in the Dirichlet--Neumann case, the higher fiber eigenvalues do not contribute to the spectrum below \(h\). In the Robin setting this follows from an exact identity, which also identifies the threshold $2mh=1-2\beta(h)$.

Throughout, $\|\cdot\|$ denotes the norm of $L^2((R,1),r\,dr)$.

\begin{proposition}[Exact identity]\label{ExactIdentity}
Let $w=u_m f$ with $f\in H^1(R,1)$ and $f(R)=0$. Then
\begin{equation}\label{ExactIdentityEq}
q^{\mathrm{Rob}}_{m,h}(w)
= h\|w\|^2
+h^2\int_R^1 u_m^2\,|f'|^2\,r\,dr
+\frac h2\bigl(2mh-1+2\beta(h)\bigr)|w(1)|^2 .
\end{equation}
\end{proposition}

\begin{proof}
Write $u=u_m$. Since $w'=u'f+uf'$, we have $|w'|^2=u^2|f'|^2+u'\,(u|f|^2)'$, hence
\[
h^2\int_R^1|w'|^2\,r\,dr
=h^2\int_R^1u^2|f'|^2\,r\,dr+h^2\int_R^1u'\,(u|f|^2)'\,r\,dr .
\]
Integrating the last term by parts and using $f(R)=0$,
\[
h^2\int_R^1u'\,(u|f|^2)'\,r\,dr
=h^2u'(1)u(1)|f(1)|^2-\int_R^1 h^2\frac{(ru')'}{r}\,u|f|^2\,r\,dr .
\]
By \eqref{um-eq}, $-h^2(ru')'/r=(h-V_m)u$, so the last integral equals $\int_R^1(h-V_m)|w|^2\,r\,dr$. Adding $\int_R^1V_m|w|^2\,r\,dr$ and the Robin term $\beta(h)h|w(1)|^2$, and using
\[
h^2u'(1)u(1)|f(1)|^2=h^2\frac{u'(1)}{u(1)}|w(1)|^2=\frac h2(2mh-1)|w(1)|^2,
\]
we obtain \eqref{ExactIdentityEq}.
\end{proof}

\begin{proposition}[Reduction to the first spectral band]\label{ReductionFirstBand}
For every $h>0$, every $\beta(h)\in\mathbb R$ and every $m\in\mathbb Z$:
\begin{enumerate}
\item[\textnormal{(i)}] if $2mh\ge 1-2\beta(h)$, then $\lambda_0(H_{m,h})>h$;
\item[\textnormal{(ii)}] $\lambda_j(H_{m,h})>h$ for every $j\ge1$.
\end{enumerate}
\end{proposition}

\begin{proof}
Every $w\in\mathcal D(q^{\mathrm{Rob}}_{m,h})$ can be written $w=u_mf$ with $f\in H^1((R,1))$ and $f(R)=0$, so \eqref{ExactIdentityEq} applies.

\smallskip
\noindent
\textbf{(i)} If $2mh\ge1-2\beta(h)$, the last term in \eqref{ExactIdentityEq} is nonnegative, so
\[
q^{\mathrm{Rob}}_{m,h}(w)\ge h\|w\|^2+h^2\int_R^1u_m^2|f'|^2\,r\,dr .
\]
Equality $q^{\mathrm{Rob}}_{m,h}(w)=h\|w\|^2$ would force $f'\equiv0$, hence $f\equiv0$ since $f(R)=0$. As $H_{m,h}$ has compact resolvent, $\lambda_0(H_{m,h})$ is attained by an eigenfunction, so $\lambda_0(H_{m,h})>h$.

\smallskip
\noindent
\textbf{(ii)} Let
\[
F=\{w\in\mathcal D(q^{\mathrm{Rob}}_{m,h}):\ w(1)=0\},
\]
a subspace of codimension one. For $w\in F$ the boundary term in \eqref{ExactIdentityEq} vanishes, so
\[
q^{\mathrm{Rob}}_{m,h}(w)=h\|w\|^2+h^2\int_R^1u_m^2|f'|^2\,r\,dr .
\]
The restriction of $q^{\mathrm{Rob}}_{m,h}$ to $F$ is the form of the Dirichlet--Dirichlet realization of $H_{m,h}$ on $(R,1)$, which has compact resolvent. Its lowest eigenvalue $\lambda_0^{DD}$ is attained, and by the argument in (i) it satisfies $\lambda_0^{DD}>h$. By the min--max principle,
\[
\lambda_1(H_{m,h})\ge\lambda_0^{DD}>h,
\]
and therefore $\lambda_j(H_{m,h})>h$ for all $j\ge1$.
\end{proof}

Consequently, only the first band function \(\lambda_0(H_{m,h})\) contributes to the counting function. It follows that
\begin{equation}\label{NB-OF-EIG2}
N(h)
=\#\Bigl\{
m\in\mathbb Z:
\lambda_0(H_{m,h})<h
\Bigr\}.
\end{equation}
By Proposition~\ref{ReductionFirstBand}(i), only integers $m$ with $2mh<1-2\beta(h)$ can contribute. The remainder of the paper is devoted to identifying, among these, the angular momenta for which $\lambda_0(H_{m,h})<h$ and deriving the corresponding asymptotic formula for $N(h)$.

When $1-2\beta-R_e^2>0$, let $\varepsilon>0$ such that $\varepsilon<\frac12 \min\Big\{ 1-2\beta-R_e^2,\,R_e^2-R^2\Big\}$ for $h$ sufficiently small. We define the \textbf{Robin-favorable interval} by
\begin{equation}\label{FavInt}
I^{\mathrm{Rob}}_{\varepsilon,h}
=
\left[
\frac{R_e^2+\varepsilon}{2h},
\frac{1-2\beta(h)-\varepsilon}{2h}
\right]
\end{equation}
The analysis of the fiber operators shows that the spectral problem reduces to identifying those angular momenta $m$ for which the first eigenvalue of $H_{m,h}$ lies below $h$. This will be achieved by distinguishing favorable and unfavorable regions.

\section{Reduction to favorable and unfavorable angular momenta}\label{DecFavUnfav}
\subsection{Favorable Angular Momenta}\label{ssec:fav}

We construct a quasi-mode localized near the outer boundary and show that its Rayleigh quotient is strictly below $h$, implying $\lambda_0(H_{m,h}) < h.$
\begin{proposition}\label{FavAngMom}Let $\beta\in\R$, $0<R<1$ and $0<\varepsilon<\frac12(1-2\beta-R_e^2)$, there exists $h_0$ such that if $h\in(0,h_0]$ and $R_e^2+\varepsilon<2mh<1-2\beta-\varepsilon$ then
\[
\lambda_0(H_{m,h})<h.
\]
\end{proposition}

\begin{proof}
Note that $r_\ast$ satisfies 
\begin{equation}\label{ineq-fav}
R^2+\varepsilon<R_e^2+\varepsilon< r_\ast^2<1-2\beta-\varepsilon.
\end{equation}
For $\beta<0$, condition \eqref{ineq-fav} reads $r_\ast^2<1+2|\beta|-\varepsilon$, so it allows $r_\ast>1$. The argument below does not use $r_\ast<1$, and in that case $u_m$ is localized at the outer boundary.

The function $\phi$ defined in \eqref{Defphi} satisfies $\phi(1)-\phi(R)=(1-R^2)-2r^2_\ast |\ln R|$. Since $r_\ast^2>R_e^2 +\varepsilon$, it follows that
\begin{equation}\label{eqt1onphi}
\phi(1)<\phi(R)-2\varepsilon |\ln R|.
\end{equation}
Choose $\eta$ such that $R+\eta<r_\ast$ and
\begin{equation}\label{eqt2onphi}
\phi(1)<\phi(R+\eta)-\varepsilon |\ln R|.
\end{equation}
This is possible due to $\eqref{eqt1onphi}$ and the uniform continuity of $\phi$. 

Note that $\phi$ is decreasing on $(0,r_\ast]$. Since $R+\eta<r_\ast$, we have $[R,R+\eta]\subset (0,r_\ast]$, hence we obtain 
\begin{equation}\label{phidecreasing}
\phi(R+\eta)\leq \phi(r)\leq \phi(R).
\end{equation}

Let $\chi$ be a smooth cut-off such that
\[
0\leq \chi\leq 1,\qquad \chi = 1 \ \text{on } [R+\eta,1], \qquad \chi=0 \text{ near } R,\qquad |\chi'|\le \frac{C}{\eta},
\]
and set $f=\chi u_m$, with $u_m$ as in \eqref{u}. Since $H_{m,h}u_m=hu_m$ by \eqref{um-eq},
$f(R)=0$ and $\chi(1)=1$, an integration by parts gives
\[
q^{\mathrm{Rob}}_{m,h}(f)=h\|f\|^2+h^2\|\chi'u_m\|^2+h\,u_m(1)^2\Bigl(h\frac{u_m'(1)}{u_m(1)}+\beta(h)\Bigr).
\]
Since $h\,u_m'(1)/u_m(1)=mh-\frac12=\frac{r_\ast^2-1}{2}$ and $u_m(1)^2=e^{-\phi(1)/(2h)}$, this reads
\begin{equation}\label{QF-h}
q^{\mathrm{Rob}}_{m,h}(f)-h\|f\|^2=h^2\|\chi'u_m\|^2+\frac h2\bigl(r_\ast^2-1+2\beta(h)\bigr)e^{-\frac{\phi(1)}{2h}} .
\end{equation}
By \eqref{eqt2onphi} and \eqref{phidecreasing},
\begin{equation}\label{estimate1}
\|\chi'u_m\|^2\le\frac{C}{\eta^2}\int_R^{R+\eta}e^{-\frac{\phi(r)}{2h}}\,r\,dr
\le\frac{C}{\eta}e^{-\frac{\phi(R+\eta)}{2h}}
\le\frac{C}{\eta}e^{-\frac{\phi(1)+\varepsilon|\ln R|}{2h}} .
\end{equation}
Moreover $r_\ast^2<1-2\beta-\varepsilon$, so $r_\ast^2-1+2\beta(h)\le-\varepsilon/2$ for $h$ small. Hence
\begin{equation}\label{equation-inv-rstar}
q^{\mathrm{Rob}}_{m,h}(f)-h\|f\|^2\le h\,e^{-\frac{\phi(1)}{2h}}
\Bigl(\frac{Ch}{\eta}e^{-\frac{\varepsilon|\ln R|}{2h}}-\frac{\varepsilon}{4}\Bigr)<0
\end{equation}
for $h$ small. By the min--max principle, $\lambda_0(H_{m,h})<h$.
\end{proof}

\subsection{Unfavorable angular momentum}\label{ssec:unfav}
In this section, we analyze the regime of unfavorable angular momenta, i.e., those values of $m$ for which the localization point $r_\ast$ defined in \eqref{rstar} lies away from the outer boundary. In this situation, the corresponding eigenfunctions are no longer localized near $r=1$, and the influence of the Dirichlet boundary at $r=R$ becomes dominant. 

Our goal is to show that, in this regime, the ground state energy of the fiber operator satisfies
\[
\lambda_0(H_{m,h}) > h,
\]
so that these angular momenta do not contribute to the counting function.

Throughout this subsection we assume $\beta(h)\ge 0$; the case $\beta<0$ is treated in Section~\ref{sec:neg}. Let $H^{\mathrm{DN}}_{m,h}$
be the operator on $L^2((R,1),r\,dr)$ associated with the form
\begin{equation}\label{eq:DN}
q^{\mathrm{DN}}_{m,h}(u)=\int_R^1\bigl(h^2|u'|^2+V_m|u|^2\bigr)\,r\,dr,
\qquad
u\in D(q^{\mathrm{DN}}_{m,h})=\{u\in H^1((R,1)) : u(R)=0\},
\end{equation}
that is, the fiber operator of the Dirichlet--Neumann problem (Dirichlet
condition at $r=R$, Neumann condition $u'(1)=0$). Since $\beta(h)\ge0$, we have
\begin{equation}\label{eq:compare}
q^{\mathrm{Rob}}_{m,h}(u)=q^{\mathrm{DN}}_{m,h}(u)+\beta(h)\,h\,|u(1)|^2
\;\ge\; q^{\mathrm{DN}}_{m,h}(u)
\end{equation}
on the common form domain, so by the min--max principle
\begin{equation}\label{eq:minmax}
\lambda_0(H_{m,h})\;\ge\;\lambda_0(H^{\mathrm{DN}}_{m,h}).
\end{equation}

\begin{proposition}\label{Unfavorable}
Let $\varepsilon \in \left(0,\frac{R_e^2-R^2}{2}\right)$.
There exists $h_0>0$ such that, for all $h\in(0,h_0]$
and all $m\in\mathbb{Z}$ with
\[
    R^2+\varepsilon < 2mh < R_e^2-\varepsilon,
\]
we have
\[
    \lambda_0(H_{m,h})>h.
\]
\end{proposition}

 To prove Proposition \ref{Unfavorable}, we need an estimate about the decay of the ground state.
\begin{proposition}\label{Agm}
Assume $\beta(h)\ge0$. Given $0<R<1$, $0<\varepsilon<\frac{1-R}{2}$ and $\gamma \in (0,1)$, there exist constants $C_0,h_0$ such that the following holds. If $h\in(0,h_0]$, $r_\ast\in [R+\varepsilon, 1-\varepsilon]$, $\lambda_0:=\lambda_0(H_{m,h})\le h$ and $\psi$ is a normalized ground state of $H_{m,h}$, then
\[
\int_R^1 e^{\gamma \frac{\phi(r)-\phi(r_\ast)}{2h}}\,|\psi|^2\, r\,dr\leq C_0 ,
\]
where $\phi$ is the function defined in \eqref{Defphi}.
\end{proposition}
\begin{proof}
Let $\Phi(r):= \gamma \frac{\phi(r)-\phi(r_\ast)}{4h}$, so that $h^2|\Phi'|^2=\gamma^2 V_m$, since $V_m=\phi'^2/16$. Recall the quadratic form \eqref{eq:fiberform}:
\begin{equation}\label{QF-new}
q^{\mathrm{Rob}}_{m,h}(e^\Phi \psi)= \int_{R}^1 \bigl(h^2 |(e^\Phi \psi)^\prime|^2 +V_m|e^\Phi \psi|^2\bigr)r\,dr+h\beta(h) e^{2\Phi(1)} |\psi(1)|^2 .
\end{equation}
Integrating $\langle H_{m,h} \psi, e^{2\Phi} \psi\rangle=\lambda_0 \|e^{\Phi} \psi\|^2$ by parts and using that $h\psi'(1)+\beta(h)\psi(1)=0$, we obtain
\begin{equation}\label{IBP}
q^{\mathrm{Rob}}_{m,h}(e^\Phi \psi)=\lambda_0 \|e^\Phi \psi\|^2 + h^2\|\Phi' e^\Phi \psi\|^2
=\lambda_0 \|e^\Phi \psi\|^2 + \gamma^2\int_{R}^1 V_{m} |e^\Phi \psi|^2\,r\,dr.
\end{equation}
Comparing \eqref{QF-new} and \eqref{IBP}, with $w=e^\Phi \psi$,
\[
h^2\|w'\|^2 +\int_{R}^1\bigl((1-\gamma^2) V_m -\lambda_0\bigr) |w|^2\, r\,dr=-h\beta(h)\, e^{2\Phi(1)}|\psi(1)|^2\le 0,
\]
since $\beta(h)\ge0$. Let $C_1^2=32/(1-\gamma^2)$ and $A=\{r\in(R,1):|\phi'(r)|\le C_1\sqrt h\}$. Since
$\phi'(r)=2(r-r_\ast)(r+r_\ast)/r$ and $r_\ast\in[R+\varepsilon,1-\varepsilon]$, we have $A\subset\{|r-r_\ast|\le C_2\sqrt h\}$
for some $C_2>0$. By Taylor's formula, $|\phi(r)-\phi(r_\ast)|\le C_3h$ on $A$, so $|\Phi|\le C_4$ on $A$
and $\int_A|w|^2r\,dr\le e^{2C_4}$. On $(R,1)\setminus A$,
$(1-\gamma^2)V_m-\lambda_0=(1-\gamma^2)\frac{\phi'^2}{16}-\lambda_0\ge 2h-h=h$, while on $A$ we have
$(1-\gamma^2)V_m-\lambda_0\ge-h$. Therefore
\[
h\int_{(R,1)\setminus A}|w|^2\,r\,dr\le h\int_A|w|^2\,r\,dr\le h\,e^{2C_4},
\]
and $\int_R^1|w|^2r\,dr\le 2e^{2C_4}=:C_0$.
\end{proof}
 Let $H_{m,h}^{\rm Dir}$ denote the operator $H_{m,h}$ on $(R,+\infty)$ with Dirichlet boundary conditions at $r=R$.
  The following estimate is the Dirichlet ingredient of \cite{FournaisKachmar}; we recall its proof for the reader's convenience.
\begin{proposition}\label{Lambda-Dir}
Given $R, \varepsilon>0$ there exist constants $c_\ast$, $c^\ast$ and $h_0$ such that if $m\in \mathbb{N}\!\!$, $h\in (0,h_0]$, $R+\varepsilon< r_*< 1-\varepsilon$
 then we have 
\[
h+c_\ast Rh^{1/2}e^{-\frac{\phi(R)-\phi(r_\ast)}{2h}}\leq \lambda_0 (H_{m,h}^{\rm Dir}) \leq h+ c^\ast R h^{1/2}e^{-\frac{\phi(R)-\phi(r_\ast)}{2h}}
\]
where $\lambda_0 (H_{m,h}^{\rm Dir})$ is the lowest eigenvalue of $H_{m,h}^{\rm Dir}$.
\end{proposition}
\begin{proof}
\emph{Scaling.} The change of variables $r = Rs$ shows that $H^{\rm Dir}_{m,h}$
on $(R,\infty)$ is unitarily equivalent to $R^2 H^{\rm Dir}_{m,\tilde{h}}$ on
$(1,\infty)$ with $\tilde{h}= h/R^2$. Under this map
$\widetilde{r_*}= r_*/R$ and
$\frac{\phi(R)-\phi(r_*)}{2h} = \frac{\tilde{\phi}(1)-\tilde{\phi}(\widetilde{r_*})}{2\tilde{h}}$,
where $\tilde{\phi}$ denotes ${\phi}$ built from $\widetilde{r_*}$, and
$R^2 (\tilde{h})^{1/2} = R h^{1/2}$. It therefore suffices to treat
$R=1$, $1+\varepsilon < r_* \le C$, in which case $\phi(1)=1$.

\emph{Test function.} Let $\rho_1 = \frac{1+r_*}{2}$ and choose
$\chi\in C^\infty(\mathbb R;[0,1])$ with $\chi = 1$ on $(-\infty,\rho_1]$,
$\chi=0$ on $[r_*,\infty)$, and $|\chi'|\le C/\varepsilon$. Set
\[
  f(r) = \int_1^r \chi(\rho)\,\rho^{-1-2m} e^{\rho^2/(2h)}\,d\rho
       = \int_1^r \chi(\rho)\,\rho^{-1} e^{\phi(\rho)/(2h)}\,d\rho,
  \qquad v = f\,u_m .
\]
Then $v(1)=0$, and $v$ decays at infinity. Since $\phi'(1) = 2-2r_*^2 \le -c<0$ and $\phi$ is decreasing
on $[1,r_*]$,
\[
  c_0 h\, e^{\frac{1}{2h}} \le f(r) \le C h\, e^{\frac{1}{2h}}
  \quad (r \ge 1+2h).
\]

\emph{Estimates.} Let $A = H^{\rm Dir}_{m,h}-h$. Then
$Av = -h^2\chi'(r)\, r^{-1-m} e^{r^2/(4h)}$, and an integration by parts
gives
\[
  \langle Av, v\rangle = h^2\int_1^\infty \chi^2(\rho)\,\rho^{-1}
  e^{\phi(\rho)/(2h)}\,d\rho .
\]
By Laplace's method at the endpoint $\rho = 1$,
\[
  c_1 h^3 e^{\frac{1}{2h}} \le \langle Av,v\rangle \le C_1 h^3 e^{\frac1{2h}} .
\]
Since $\operatorname{supp}\chi' \subset [\rho_1, r_*]$, where
$\phi \le \phi(\rho_1) < 1$,
\[
  \|Av\|^2 \le C h^4 e^{\phi(\rho_1)/(2h)} .
\]
By the bounds on $f$ and \eqref{LaplaceMethod},
\[
  \|v\|^2 \asymp h^2 e^{1/h} \int_1^\infty |u_m|^2 r\,dr
  \asymp h^{5/2} e^{\frac{2-\phi(r_*)}{2h}} .
\]
Hence, with $\tau = \langle Av,v\rangle/\|v\|^2$ and
$\sigma^2 = \|Av\|^2/\|v\|^2 - \tau^2$,
\[
  \tau \asymp h^{1/2} e^{-\frac{1-\phi(r_*)}{2h}}, \qquad
  \sigma^2 \le C h^{3/2} e^{-\frac{1-\phi(r_*)}{2h}}
             e^{-\frac{1-\phi(\rho_1)}{2h}} .
\]

\emph{Conclusion.} Write $E_j:=\lambda_j(H^{\rm Dir}_{m,h})-h$, $j=0,1$.
By the harmonic approximation,
$E_1 = \lambda_1^{\rm Dir} - h = 2h + O(h^{3/2})$ (see \cite{FournaisHelfferBook}), so
$E_1 - \tau \ge h$ for $h$ small. Temple's inequality gives
\[
  E_0 \ge \tau - \frac{\sigma^2}{E_1-\tau}
      \ge \tau - C h^{1/2} e^{-\frac{1-\phi(r_*)}{2h}}
          e^{-\frac{1-\phi(\rho_1)}{2h}}
      \ge c_* h^{1/2} e^{-\frac{1-\phi(r_*)}{2h}},
\]
and the min--max principle gives $E_0 \le \tau \le
c^* h^{1/2} e^{-\frac{1-\phi(r_*)}{2h}}$. Since
$1-\phi(r_*) = \phi(1)-\phi(r_*)$, scaling back proves the claim.
\end{proof}

\begin{proof}[Proof of Proposition \ref{Unfavorable}]
The condition on $m$ reads as 
\[
R^2+\varepsilon<r_\ast^2<R_e^2-\varepsilon <1-\varepsilon 
\]
Hence $r_* \in [R+\varepsilon', 1-\varepsilon']$ with $\varepsilon' := \varepsilon/2$, so Propositions~\ref{Agm} and~\ref{Lambda-Dir} apply with $\varepsilon'$ in place of $\varepsilon$.

In a similar manner to the proof of \eqref{eqt1onphi}, using that $r_\ast^2<R_e^2-\varepsilon$, we get 
\begin{equation}\label{LBphi}
\phi(1)>\phi(R)+2\varepsilon |\ln R|.
\end{equation}
Since \(\phi\) is continuous on the compact interval
\([R,1]\), it is uniformly continuous.
Hence, choosing \(\eta\) sufficiently small, such that $1-\eta >r_\ast$, we have
\[
|\phi(1-\eta)-\phi(1)|
<
\varepsilon |\ln R|.
\]
This together with \eqref{LBphi} gives
\begin{equation}\label{est9}
\phi(1-\eta)>\phi(R)+\varepsilon |\ln R|.
\end{equation}

By \eqref{est9}, $\phi(1-\eta) > \phi(R)$, so we may choose $\gamma \in (0,1)$ with
$$\gamma\big(\phi(1-\eta)-\phi(r_*)\big) > \phi(R)-\phi(r_*).$$
Since $r_*$ ranges over a compact set and \eqref{est9} holds with the margin $\varepsilon|\ln R|$, $\eta$ and $\gamma$ can be chosen depending only on $R$ and $\varepsilon$, not on $m$.
Assume, by contradiction, that $\lambda_0(H_{m,h}) \le h$; then Proposition~\ref{Agm} applies with this $\gamma$.

Let $\psi$ be a normalized ground state of $H_{m,h}$. We introduce $f=\chi \psi$ where $\chi$ satisfies
\[
0\leq \chi\leq 1,\qquad {\rm supp}\chi \subset (-\infty,1-\frac\eta2), \qquad \chi \Big|_{(-\infty,1-\eta)}=1
\]
 Recall the quadratic form $q^{\mathrm{Rob}}_{m,h}$ defined in \eqref{eq:fiberform}. An integration by parts yields
\begin{equation}\label{est7}
q^{\mathrm{Rob}}_{m,h}(f)=\lambda_0 (H_{m,h})\|f\|^2+h^2\|\chi^\prime \psi\|^2.
\end{equation}
Note that the boundary term vanishes since $\psi(R)=0$ and $\chi(1)=0.$
Since $0\le\chi\le 1$, $\chi=1$ on $(R,1-\eta)$, and $\|\psi\|=1$, we have
$\|f\|^2\le 1$. A lower bound is obtained as follows. Since $1-\eta>r_*$,
$\phi$ is increasing on $[1-\eta,1]$, so $\phi(r)\ge\phi(1-\eta)$ there. Hence,
by Proposition~\ref{Agm},
\begin{equation}\label{47new}
\begin{aligned}
1-\|f\|^2
  &=\int_{1-\eta}^{1}\bigl(1-\chi^2\bigr)|\psi|^2\,r\,dr
   \;\le\;\int_{1-\eta}^{1}|\psi|^2\,r\,dr\\
  &\le e^{-\gamma\frac{\phi(1-\eta)-\phi(r_*)}{2h}}
   \int_{1-\eta}^{1} e^{\gamma\frac{\phi(r)-\phi(r_*)}{2h}}|\psi|^2\,r\,dr
  \;\le\; C_0'\,e^{-\gamma\frac{\phi(1-\eta)-\phi(r_*)}{2h}} .
\end{aligned}
\end{equation}
Set $\varepsilon_h:=e^{-\gamma\frac{\phi(1-\eta)-\phi(r_*)}{2h}}$. Then
\begin{equation}\label{47norm}
  1-C_0'\varepsilon_h\;\le\;\|f\|^2\;\le\;1 ,
\end{equation}
and in the same way
\begin{equation}\label{47chi}
  \|\chi'\psi\|^2\le\frac{C}{\eta^2}\int_{1-\eta}^{1}|\psi|^2\,r\,dr
  =O(\varepsilon_h).
\end{equation}

Since $f(R)=f(1)=0$ (so the Robin term vanishes), $f$ is an admissible function for
$H^{\mathrm{Dir}}_{m,h}$, and Proposition \ref{Lambda-Dir} gives
\begin{equation}\label{48new}
  q^{\mathrm{Rob}}_{m,h}(f)\;\ge\;\lambda_0(H^{\mathrm{Dir}}_{m,h})\,\|f\|^2 .
\end{equation}
By the IMS-type identity \eqref{est7},
$q^{\mathrm{Rob}}_{m,h}(f)=\lambda_0(H_{m,h})\|f\|^2+h^2\|\chi'\psi\|^2$.
Combining this with \eqref{48new} and dividing by $\|f\|^2\ge\tfrac12$
(valid for $h$ small, by \eqref{47norm}), we obtain
\begin{equation}\label{50new}
\begin{aligned}
  \lambda_0(H_{m,h})
  &\;\ge\;\lambda_0(H^{\mathrm{Dir}}_{m,h})-\frac{h^2\|\chi'\psi\|^2}{\|f\|^2}\\
  &\;\ge\; h+c_*R\,h^{1/2}e^{-\frac{\phi(R)-\phi(r_*)}{2h}}-O(\varepsilon_h).
\end{aligned}
\end{equation}
By the choice of $\gamma$, we have $O(\varepsilon_h)=o\!\bigl(h^{1/2}e^{-\frac{\phi(R)-\phi(r_*)}{2h}}\bigr)$.
Hence, for $h$ small, \eqref{50new} yields $\lambda_0(H_{m,h})>h$.

\end{proof}

\section{Proof of the Main Theorem for \texorpdfstring{$\beta\ge0$}{beta>=0}}\label{sec:proofmain}

Throughout this section we assume $\beta(h)\ge0$ for small $h$, which holds when $\beta>0$. The case $\beta = 0$ with $\beta(h)$ negative or changing sign is treated in Remark~\ref{rem:beta0}. We set
\[
x(h)=1-R_e^2-2\beta(h),\qquad x=\lim_{h\to0}x(h)=1-R_e^2-2\beta .
\]
We prove $2h\,N(h)\to x_+$.
\begin{proof}
By Proposition~\ref{ReductionFirstBand},
\[
N(h)=\#\{m\in\mathbb Z:\ \lambda_0(H_{m,h})<h\}.
\]

\medskip
\noindent\textbf{Lower bound.}
If $x\le 0$ there is nothing to prove. Let $x>0$ and fix $\varepsilon$ as in \eqref{FavInt}. For $h$ small, $|\beta(h)-\beta|<\varepsilon/4$, so every $m\in I^{\rm Rob}_{\varepsilon,h}$ satisfies $R_e^2+\varepsilon/2<2mh<1-2\beta-\varepsilon/2$, and Proposition~\ref{FavAngMom} (with $\varepsilon/2$ in place of $\varepsilon$) gives $\lambda_0(H_{m,h})<h$.
Since an interval of length $L$ contains at least $L-1$ integers, we get
\[
N(h)\ \ge\ \bigl|I^{\mathrm{Rob}}_{\varepsilon,h}\bigr|-1
=\frac{x(h)-2\varepsilon}{2h}-1 .
\]

\medskip
\noindent\textbf{Upper bound.}
Fix $0<\varepsilon<\tfrac12(R_e^2-R^2)$ and let $m\in\mathbb Z$ satisfy $\lambda_0(H_{m,h})<h$. We show that $2mh$ lies in a small set.

\smallskip
\emph{(a) Large $m$.} By Proposition~\ref{ReductionFirstBand}(i), $2mh<1-2\beta(h)$.

\smallskip
\emph{(b) Small $m$.} Since $\beta(h)\ge0$, the Robin term in \eqref{eq:fiberform} is nonnegative, so
$q^{\mathrm{Rob}}_{m,h}(u)\ge\int_R^1\bigl(h^2|u'|^2+V_m|u|^2\bigr)r\,dr$.
If $2mh\le R^2-4R\sqrt h$, then $\frac{mh}{r}-\frac r2\le-2\sqrt h$ on $(R,1)$ for small $h$, hence $V_m\ge4h>h$ on $(R,1)$. Indeed, for $m\ge0$ the supremum of $\frac{mh}{r}-\frac r2$ over $(R,1)$ is attained at $r=R$ and equals $\frac{mh}{R}-\frac R2\le-2\sqrt h$; for $m<0$ it is at most $-R/2\le-2\sqrt h$. Then $q^{\mathrm{Rob}}_{m,h}(u)\ge4h\|u\|^2$, so $\lambda_0(H_{m,h})>h$. Hence $2mh>R^2-4R\sqrt h$.

\smallskip
\emph{(c) Intermediate $m$.} Let $R^2+\varepsilon<2mh<R_e^2-\varepsilon$. Then $\lambda_0(H_{m,h})>h$ by Proposition~\ref{Unfavorable}.

\smallskip
Combining (a)--(c), $2mh$ belongs to
\[
\bigl(R^2-4R\sqrt h,\ R^2+\varepsilon\bigr]\ \cup\ \bigl[R_e^2-\varepsilon,\ 1-2\beta(h)\bigr).
\]
An interval of length $L$ contains at most $L+1$ integers, so
\[
N(h)\le\frac{\varepsilon+4R\sqrt h}{2h}+1+\frac{\bigl(x(h)+\varepsilon\bigr)_+}{2h}+1
\le\frac{x(h)_++2\varepsilon}{2h}+\frac{2R}{\sqrt h}+2 .
\]

\medskip
\noindent\textbf{Conclusion.}
Multiplying by $2h$ and using $x(h)\to x$,
\[
(x-2\varepsilon)_+\ \le\ \liminf_{h\to0}2hN(h)\ \le\ \limsup_{h\to0}2hN(h)\ \le\ x_++2\varepsilon .
\]
Letting $\varepsilon\to0$ gives $2hN(h)\to x_+=\ell(\beta)$, since the boundary correction $-4\beta$ is non-positive for $\beta\geq 0$.
\end{proof}

\begin{remark}
The assumption $\beta(h)\ge0$ is used only in the upper bound, steps (b) and (c), where the Robin term is nonnegative (in step (c) through Proposition~\ref{Agm}). For $\beta<0$ the boundary term is negative and these two steps require a new argument; this is done in Section~\ref{sec:neg}, where Theorem~\ref{thm:main_Robin} is proved for $\beta<0$.
\end{remark}

\section{The attractive case \texorpdfstring{$\beta<0$}{beta<0}}\label{sec:neg}

Throughout this section $\beta(h)\to\beta<0$, and we write
\[
  \bh:=|\beta(h)|=-\beta(h)\longrightarrow|\beta|>0 ,
  \qquad
  \beta_c=\frac{1-R_e^2}{2},
  \qquad
  \|\cdot\|=\|\cdot\|_{L^2((R,1),r\dd r)} .
\]
For $h$ small we may and do assume $\tfrac12|\beta|\le \bh\le 2|\beta|$.
The fiber form \eqref{eq:fiberform} reads
\begin{equation}\label{eq:form-neg}
  \q(u)=\int_R^1\bigl(h^2|u'|^2+V_m|u|^2\bigr)r\dd r-\bh\,h\,|u(1)|^2 ,
  \qquad
  V_m(r)=\Bigl(\frac{mh}{r}-\frac r2\Bigr)^2 ,
\end{equation}
so the boundary term is now \emph{negative} (attractive).

The aim of this section is to prove Theorem~\ref{thm:main_Robin} for $\beta<0$. The function $\ell$ is continuous and piecewise linear, with slope $-4$ on
$(-\infty,-\beta_c)$, slope $-2$ on $(-\beta_c,\beta_c)$ and $\ell\equiv0$ on
$[\beta_c,\infty)$; see Figure~\ref{fig:density-vs-beta}. 

\paragraph{Reduction to the ground state.}
Since $\lambda_0(H_{m,h})$ can be as low as $-\beta^2$
(Corollary~\ref{cor:bottom}), one might expect the attractive boundary
to push several fiber eigenvalues below $h$. This does not happen:
Proposition~\ref{ReductionFirstBand} was stated for every $\beta(h)\in\mathbb{R}$,
and the proof of part~(ii) only uses the form on
$F=\{w : w(1)=0\}$, where the Robin term $-b_h h|w(1)|^2$ vanishes.
Hence, also for $\beta<0$,
\[
  \lambda_j(H_{m,h}) \ge \lambda_0(H^{DD}_{m,h}) > h
  \qquad (j\ge 1),
\]
so each fiber contributes at most one eigenvalue below $h$, and \eqref{NB-OF-EIG2} holds. Proposition~\ref{ReductionFirstBand}(i) carries over as well, but it says
nothing on the range $2mh<1+2b_h$, where $\lambda_0$ may lie at
distance $O(1)$ below $h$ (Proposition~\ref{prop:bs} below).
For $\beta<0$, therefore, all the work goes into $\lambda_0$.

\subsection{A sharp semiclassical trace inequality}

\begin{lemma}[Sharp trace inequality]\label{lem:trace}
Let $b>0$, $h>0$, $\delta\in(0,\frac12]$ and $u\in H^1((1-\delta,1))$. Then
\begin{equation}\label{eq:trace}
  b\,h\,|u(1)|^2\le h^2\int_{1-\delta}^1|u'|^2\,r\dd r
  +\Bigl(b^2+\frac{4bh}{\delta}\Bigr)\int_{1-\delta}^1|u|^2\,r\dd r .
\end{equation}
The constant $b^2$ is optimal as $h\to0$; see the proof of Proposition~\ref{prop:bs}.
\end{lemma}

\begin{proof}
Let $\rho(r)=(r-1+\delta)/\delta$, so that $\rho(1-\delta)=0$, $\rho(1)=1$ and $0\le\rho\le1$. Then
\[
  |u(1)|^2=\int_{1-\delta}^1\frac{\dd}{\dd r}\bigl(\rho\,r\,|u|^2\bigr)\dd r
  =\int_{1-\delta}^1\Bigl(\frac r\delta+\rho\Bigr)|u|^2\dd r
   +2\int_{1-\delta}^1\rho\,r\,\mathrm{Re}(\bar u u')\dd r .
\]
By Young's inequality, $2bh\,\rho r|u||u'|\le h^2r|u'|^2+b^2\rho^2r|u|^2\le h^2r|u'|^2+b^2r|u|^2$.
Moreover $r/\delta+\rho\le 2/\delta$ and $\dd r\le 2\,r\dd r$ on $[1-\delta,1]$, since $\delta\le\frac12$.
Multiplying by $bh$ and combining gives \eqref{eq:trace}.
\end{proof}

\begin{corollary}[Splitting the form]\label{cor:split}
For $u\in D(\q)$ and $\delta\in(0,\min\{\frac12,1-R\})$,
\begin{equation}\label{eq:split}
  \q(u)\ \ge\ \int_R^{1-\delta}\bigl(h^2|u'|^2+V_m|u|^2\bigr)r\dd r
  +\int_{1-\delta}^1\Bigl(V_m-\bh^2-\frac{4\bh h}{\delta}\Bigr)|u|^2\,r\dd r .
\end{equation}
In particular $\lambda_0(H_{m,h})\ge-\bh^2-4\bh h/\delta$ for every $m\in\mathbb Z$.
\end{corollary}

\begin{proof}
Apply \eqref{eq:trace} with $b=\bh$ to $u|_{(1-\delta,1)}$ in \eqref{eq:form-neg}.
\end{proof}

Now, the goal is to study the sign of
$V_m-\bh^2$ in the layer, which the next elementary lemma settles.

\begin{lemma}[The potential at the wall]\label{lem:wall}
Let $\varepsilon>0$, $\delta\in(0,\min\{\frac14,\frac\varepsilon8\}]$ and let $m\in\mathbb Z$ satisfy
$1-2mh\ge 2\bh+\varepsilon$. Then, for $h\le\varepsilon\delta/8$,
\[
  V_m(r)-\bh^2-\frac{4\bh h}{\delta}\ \ge\ \frac{\varepsilon^2}{16}
  \qquad\text{for all }r\in[1-\delta,1].
\]
\end{lemma}

\begin{proof}
Let $g(r)=\frac r2-\frac{mh}{r}$, so that $V_m=(g(r))^2$ and $g(1)=\frac{1-2mh}2$. Note that $2mh<1$ since $1-2mh\ge 2\bh+\varepsilon>0$. So either $m\leq 0$ or $0<mh<\frac12$.
\begin{itemize}
\item If $m\le0$ then $g(1)-g(r)=\frac{1-r}{2}-|m|h(\frac1r-1)\le\frac\delta2$.
\item If $0<mh<\frac12$ then $0<g'(r)=\frac12+\frac{mh}{r^2}<2$ on $[\frac34,1]$, so $g(1)-g(r)\le2\delta$.
\end{itemize}
In both cases $g(r)\ge\bh+\frac\varepsilon2-2\delta\ge\bh+\frac\varepsilon4>0$, hence
$V_m(r)\ge\bh^2+\frac{\bh\varepsilon}{2}+\frac{\varepsilon^2}{16}$.
Finally $\frac{\bh\varepsilon}{2}\ge\frac{4\bh h}{\delta}$ when $h\le\varepsilon\delta/8$.
\end{proof}

\subsection{Boundary states}

\begin{proposition}[Boundary states]\label{prop:bs}
Let $\beta<0$ and $0<\varepsilon<|\beta|$. There exists $h_0>0$ such that, for
$h\in(0,h_0]$ and every $m\in\mathbb Z$ with $|1-2mh|\le2|\beta|-\varepsilon$,
\[
  \lambda_0(H_{m,h})\le-\frac{\varepsilon|\beta|}{8}<h .
\]
\end{proposition}

\begin{proof}
For $h$ small, $\bh\ge|\beta|-\varepsilon/4$, so $a:=|1-2mh|\le2\bh-\varepsilon/2$.
Also $|mh|\le M:=1+|\beta|$, and therefore $|V_m'|\le L$ on $[\frac12,1]$ for some $L=L(\beta)$.
Fix $\delta=\min\{\frac14,\frac{1-R}{4}\}$ and a cut-off $\chi\in C^\infty$ with
$0\le\chi\le1$, $\chi=1$ on $[1-\delta,1]$, $\chi=0$ on $[R,1-2\delta]$ and $|\chi'|\le2/\delta$.
Set
\[
  w(r)=\chi(r)\,e^{-\bh(1-r)/h}\in D(\q).
\]
All integrals below over $[1-2\delta,1-\delta]$ are $O(e^{-\bh\delta/h})$. With $t=1-r$:
\[
  \|w\|^2=\int_0^\delta e^{-2\bh t/h}(1-t)\dd t+O(e^{-\bh\delta/h})
  =\frac{h}{2\bh}\bigl(1+O(h)\bigr).
\]
Since $w'=\frac{\bh}{h}w$ on $[1-\delta,1]$, we have $h^2\|w'\|^2=\bh^2\|w\|^2+O(e^{-\bh\delta/h})$. Since we have
 $|V_m(r)-V_m(1)|\le L(1-r)$ and $\int_0^\infty t\,e^{-2\bh t/h}\dd t=h^2/(4\bh^2)$,
\[
  \int_R^1V_m|w|^2r\dd r\le V_m(1)\|w\|^2+\frac{Lh^2}{4\bh^2}
  \le\bigl(V_m(1)+Ch\bigr)\|w\|^2 .
\]
Finally $\bh h|w(1)|^2=\bh h=2\bh^2\|w\|^2(1+O(h))$. Altogether,
\[
  \frac{\q(w)}{\|w\|^2}\le\bh^2+V_m(1)-2\bh^2+Ch=V_m(1)-\bh^2+Ch .
\]
Now $V_m(1)-\bh^2=\frac14(a-2\bh)(a+2\bh)\le-\frac14\cdot\frac\varepsilon2\cdot2\bh
=-\frac{\varepsilon\bh}{4}\le-\frac{3\varepsilon|\beta|}{16}$.
Taking $h$ small enough that $Ch\le\varepsilon|\beta|/16$, the min--max principle gives the claim.
\end{proof}

\begin{corollary}[Bottom of the spectrum]\label{cor:bottom}
For $\beta<0$,
\[
  \inf\mathrm{Spec}(P_h)=-\beta(h)^2+O(h).
\]
In particular, the spectrum below the lowest Landau level reaches an
$O(1)$ distance below $h$: for the fibers of Proposition~\ref{prop:bs},
$\lambda_0(H_{m,h})$ is not close to $h$.
\end{corollary}

\begin{proof}
The lower bound is Corollary~\ref{cor:split} with $\delta$ fixed. For the upper
bound, we take $m\!=\lfloor 1/(2h)\rfloor$, so that $V_m(1)\le h^2$, and use the
quasimode above.
\end{proof}

\subsection{Small angular momenta}

\begin{proposition}\label{prop:small}
Let $\beta<0$ and $\varepsilon\in(0,1)$. There exists $h_0>0$ such that, for
$h\in(0,h_0]$ and every $m\in\mathbb Z$ with
\[
  2mh\le R^2-4R\sqrt h\qquad\text{and}\qquad 1-2mh\ge2\bh+\varepsilon,
\]
we have $\lambda_0(H_{m,h})\ge4h>h$.
\end{proposition}

\begin{proof}
Take $\delta=\min\{\frac\varepsilon8,\frac14,\frac{1-R}{2}\}$ in \eqref{eq:split}. On $(R,1)$ we have
$V_m\ge4h$, exactly as in step~(b) of Section~\ref{sec:proofmain}. On $[1-\delta,1]$,
Lemma~\ref{lem:wall} bounds the integrand below by $\varepsilon^2/16\ge4h$.
Hence $\q(u)\ge4h\|u\|^2$.
\end{proof}

\subsection{Unfavorable angular momenta: a frozen Agmon weight}

For $\delta\in(0,1-R)$ we introduce the \emph{truncated} weight
\begin{equation}\label{eq:weight}
  \Phi_\delta(r)=\frac{\gamma}{4h}\Bigl(\phi\bigl(\min\{r,1-\delta\}\bigr)-\phi(r_*)\Bigr),
  \qquad
  h^2|\Phi_\delta'|^2=
  \begin{cases}\gamma^2V_m & \text{on }(R,1-\delta),\\ 0&\text{on }(1-\delta,1),\end{cases}
\end{equation}
where $\phi$ is as in \eqref{Defphi}; recall $V_m=\phi'^2/16$.



\begin{proposition}\label{prop:agmon}
Let $\beta<0$, $\kappa\in(0,\frac{1-R}{2})$, $\varepsilon\in(0,1)$, $\gamma\in(0,1)$ and
$\delta\in\big(0,\min\{\frac{\varepsilon}{8},\frac{\kappa}{2},\frac14\}\big]$. There exist
$C_0,h_0>0$, depending on $R,\beta,\kappa,\varepsilon,\gamma,\delta$ but not on $m$, such that:
if $h\in(0,h_0]$, $m\in\mathbb Z$ satisfies
\[ r_*\in[R+\kappa,\,1-\kappa], \qquad 1-2mh\ge 2b_h+\varepsilon, \]
$\lambda_0:=\lambda_0(H_{m,h})\le h$ and $\psi$ is a normalized ground state, then
$$\int_R^1 e^{2\Phi_\delta}|\psi|^2\,r\,dr\le C_0.$$
\end{proposition}

\begin{proof}
Let $w=e^{\Phi_\delta}\psi$. By \eqref{IBP} and \eqref{eq:weight},
\[
  \int_R^{1-\delta}\!\!\bigl(h^2|w'|^2+((1-\gamma^2)V_m-\lambda_0)|w|^2\bigr)r\dd r
  +\int_{1-\delta}^1\!\!\bigl(h^2|w'|^2+(V_m-\lambda_0)|w|^2\bigr)r\dd r
  -\bh h|w(1)|^2=0 .
\]
This is where the frozen weight matters: on $(1-\delta,1)$ no $|\Phi'|^2$
appears. We apply Lemma~\ref{lem:trace} to $w$ on $(1-\delta,1)$, then
Lemma~\ref{lem:wall}, and use $\lambda_0\le h\le\varepsilon^2/32$. Dropping
$h^2|w'|^2$ on $(R,1-\delta)$, we obtain
\begin{equation}\label{eq:agmon1}
  \int_R^{1-\delta}\bigl((1-\gamma^2)V_m-\lambda_0\bigr)|w|^2r\dd r
  +\frac{\varepsilon^2}{32}\int_{1-\delta}^1|w|^2r\dd r\ \le\ 0 .
\end{equation}
Let $C_1^2=32/(1-\gamma^2)$ and $A=\{r\in(R,1-\delta)\!:\!|\phi'(r)|\le C_1\sqrt h\}$.\! Since $\phi'(r)\!=\!\frac{2(r-r_*)(r+r_*)}{r}$ and $r_*\le1-\kappa\le 1-2\delta$, we have, for $h$ small,
$A\subset\{|r-r_*|\le C_2\sqrt h\}\subset(R,1-\delta)$, with $C_2=C_2(R,\kappa,\gamma)$. Hence $|\Phi_\delta|\le C_4$
on $A$, as in the proof of Proposition~\ref{Agm}, and
$\int_A|w|^2r\dd r\le e^{2C_4}$. On $(R,1-\delta)\setminus A$ we have
$(1-\gamma^2)V_m-\lambda_0\ge(1-\gamma^2)C_1^2h/16-h=h$. On $A$ we have
$(1-\gamma^2)V_m-\lambda_0\ge-h$. Thus \eqref{eq:agmon1} gives
\[
  h\int_{(R,1-\delta)\setminus A}|w|^2r\dd r+\frac{\varepsilon^2}{32}\int_{1-\delta}^1|w|^2r\dd r
  \le h\int_A|w|^2r\dd r\le h\,e^{2C_4},
\]
and therefore $\int_R^1|w|^2r\dd r\le e^{2C_4}(2+32h/\varepsilon^2)\leq 3e^{2C_4}=:C_0$.
\end{proof}

\begin{proposition}\label{prop:unfav}
Let $\beta<0$, $\mu=\min\{R_e^2,1-2|\beta|\}$ and assume $\mu>R^2$. Let
$0<\varepsilon<\frac{\mu-R^2}{4}\!\!$. There exists $h_0>0$ such that, for
$h\in(0,h_0]$ and every $m\in\mathbb Z$ with
\[
  R^2+\varepsilon<2mh<\mu-2\varepsilon,
\]
we have $\lambda_0(H_{m,h})>h$.
\end{proposition}

\begin{proof}
Suppose $\lambda_0(H_{m,h})\le h$. Since $r_*^2>R^2+\varepsilon$ and $r_*+R<2$, we have
$r_*-R>\varepsilon/2$; since $r_*^2<\mu-2\varepsilon$,
$1-r_*\ge \frac{1-r_*^2}{2}>\frac{1-\mu}{2}+\varepsilon\ge\varepsilon$.
Hence $r_*\in[R+\kappa,1-\kappa]$ with $\kappa:=\varepsilon/2$, and $\kappa<\frac{1-R}{2}$
because $\varepsilon<\frac{\mu-R^2}{4}<\frac{1-R}{2}$. Moreover
$1-2mh>1-\mu+2\varepsilon\ge 2|\beta|+2\varepsilon\ge 2b_h+\varepsilon$ for $h$ small.
Since $r_*^2<R_e^2-2\varepsilon$, choose $\eta<\varepsilon$ and $\gamma\in(0,1)$ as in the
proof of Proposition~\ref{Unfavorable}; they depend only on $R$ and $\varepsilon$, and $1-\eta>r_*$.
Set $\delta:=\min\{\eta,\frac{\varepsilon}{8},\frac14\}$. Then
$\delta\le\min\{\frac{\varepsilon}{8},\frac{\kappa}{2},\frac14\}$, so
Proposition~\ref{prop:agmon} applies with $(\kappa,\varepsilon,\gamma,\delta)$.

The cut-off $\chi$ in \eqref{est7} vanishes at $r=1$, so the IMS identity
\eqref{est7} and the comparison \eqref{48new} with $H^{\mathrm{Dir}}_{m,h}$ do not
involve the Robin condition. Then \eqref{50new} and Proposition~\ref{Lambda-Dir} give
$\lambda_0(H_{m,h})>h$ for $h$ small, which is a contradiction.
\end{proof}

\subsection{Proof of Theorem~\ref{thm:main_Robin} for $\beta<0$} \label{proof-beta-neg}

Let $\mu=\mu(\beta):=\min\{R_e^2,1-2|\beta|\}$ and $\ell=1+2|\beta|-\mu$. Recall, by \eqref{NB-OF-EIG2}:
$$\NRob(h)=\#\{m:\lambda_0(H_{m,h})<h\}.$$

\paragraph{Lower bound.}
\begin{itemize}
\item Proposition~\ref{FavAngMom} holds for every $\beta\in\mathbb R$. Since $1-2\beta=1+2|\beta|$, it gives $\lambda_0(H_{m,h})<h$ for $2mh\in[R_e^2+\varepsilon,\,1+2|\beta|-\varepsilon]$. 

\item Proposition~\ref{prop:bs} gives $\lambda_0<h$ for
$2mh\in[1-2|\beta|+\varepsilon,1+2|\beta|-\varepsilon]$.
\end{itemize}
Their union is
$[\mu+\varepsilon,\,1+2|\beta|-\varepsilon]$. Therefore
$\NRob(h)\ge\dfrac{\ell-2\varepsilon}{2h}-1$.

\paragraph{Upper bound.}
Let $\lambda_0(H_{m,h})<h$. Then:
\begin{enumerate}
\item $2mh<1+2\bh$, by Proposition~\ref{ReductionFirstBand}(i);
\item if $2mh\le R^2-4R\sqrt h$, then by Proposition~\ref{prop:small} we have
  $1-2mh<2\bh+\varepsilon$, i.e.\ $2mh>1-2\bh-\varepsilon\ge\mu-2\varepsilon$ for $h$ small;
\item $2mh\notin(R^2+\varepsilon,\mu-2\varepsilon)$, by Proposition~\ref{prop:unfav}.
\end{enumerate}
Hence
\[
  2mh\in\bigl(R^2-4R\sqrt h,\,R^2+\varepsilon\bigr]\ \cup\ \bigl[\mu-2\varepsilon,\,1+2\bh\bigr),
\]
and
\[
  \NRob(h)\le\frac{\varepsilon+4R\sqrt h}{2h}+\frac{1+2\bh-\mu+2\varepsilon}{2h}+2 .
\]
(If $\mu\le R^2$, which can only happen when $\beta<-\beta_c$, Proposition~\ref{prop:unfav}
is not needed and 3. is void; the same bound holds.)

\paragraph{Conclusion.}
$\ell-3\varepsilon\le\liminf 2h\NRob(h)\le\limsup2h\NRob(h)\le\ell+3\varepsilon$.
Letting $\varepsilon\to0$ proves the theorem.
\qed

\begin{remark}[The case $\beta = 0$ with $\beta(h)$ of arbitrary sign]\label{rem:beta0}
\label{rem:beta0}
Section \ref{sec:proofmain} assumes $\beta(h)\ge 0$, and the above proof assumes
$b_h \ge \tfrac12|\beta| > 0$. The remaining case $\beta(h)\to 0$
with $\beta(h)<0$ (or of changing sign) is covered by the same
arguments. The lower bound $b_h \ge \tfrac12|\beta|$ is used only in
Proposition~\ref{prop:bs} and Corollary~\ref{cor:bottom} which are not needed when
$\beta = 0$. Lemmas~\ref{lem:trace} and \ref{lem:wall} and Propositions~\ref{prop:small}--\ref{prop:unfav} hold for
every $b_h \in [0,1]$, with constants independent of $b_h$.

Let $\varepsilon>0$ be small. For $h$ small with $\beta(h)<0$ we have
$2b_h<\varepsilon$, so the hypotheses
$1-2mh \ge 2b_h+\varepsilon$ in Propositions~\ref{prop:small} and \ref{prop:agmon} hold on the
ranges considered, with $\mu = R_e^2$. The upper bound of
Section~\ref{proof-beta-neg} then gives
\[
  2h\,N(h) \le 1 - R_e^2 + 2b_h + 3\varepsilon + o(1).
\]
For $h$ with $\beta(h)\ge 0$, the bound of Section~\ref{sec:proofmain} applies.
Proposition~\ref{FavAngMom} gives the lower bound in both cases. Hence
\[
  \lim_{h\to0} 2h\,N(h) = 1-R_e^2 = \ell(0),
\]
whatever the sign of $\beta(h)$, and \eqref{eq:ell-all} holds for every
$\beta\in\mathbb R$.
\end{remark}

\begin{remark}[Relation to the Robin Laplacian with large negative parameter]
In unscaled form the boundary condition is $u'(1)=\bh h^{-1}u(1)$. This is a
Robin Laplacian with large negative parameter $\tau=-\bh/h$. The energy
$-\bh^2=h^2\cdot(-\tau^2)$ is the classical leading term $-\tau^2$
\cite{HelfferKachmar,PankrashkinPopoff}, and Lemma~\ref{lem:trace} is the one-dimensional, cut-off
version of the reduction to a boundary layer used in
\cite{KKR16,HKR17}. Weyl-type counting for the magnetic Laplacian with
de~Gennes/Robin condition in the edge-layer scaling is in
\cite{Ka08}; see also \cite{Ka06} for the ground state energy. The present regime $\beta(h)\to\beta$ is different: there the
attraction is $O(1)$ in energy, not $O(h)$.
\end{remark}

\section{Perspectives}\label{sec:perspectives}
A natural continuation of the present work is to consider Robin boundary conditions depending on the angular variable. More precisely, one may replace the constant parameter \(\beta\) by a smooth function \(\beta(\theta)\) on the outer boundary of the annulus. In this setting, the rotational symmetry is lost and the fiber decomposition is no longer available, making the analysis substantially more delicate.

Nevertheless, the asymptotic formula obtained in this paper suggests that the counting function should still be governed by the local strength of the Robin interaction. We expect, on the basis of a formal local-freezing heuristic in $\theta$, that the counting function should satisfy
\[N(h)
\sim
\frac{1}{4\pi h}
\int_{0}^{2\pi}
\ell\bigl(\beta(\theta)\bigr)\,d\theta,
\qquad h\to0,
\]
with $\ell$ as in \eqref{eq:ell-all},
though establishing this rigorously would require replacing the fiber decomposition with a more delicate microlocal analysis, since rotational symmetry is lost.
Such a result would provide a semiclassical description of the influence of non-uniform boundary interactions on the distribution of eigenvalues below the lowest Landau level.
Another natural extension concerns the strip geometry endowed with Robin boundary conditions. The techniques developed in the present work suggest that a similar asymptotic behavior should hold in that setting. This question will be addressed in future work.

\section*{Acknowledgements}
The authors are grateful to Ayman Kachmar for many helpful discussions, and for valuable
comments on earlier versions of this manuscript.

\section*{Declarations}

\noindent\textbf{Data availability.} Data sharing is not applicable to this
article, as no datasets were generated or analysed during the current study.

\noindent\textbf{Competing interests.} The authors have no competing interests
to declare that are relevant to the content of this article.

\noindent\textbf{Funding.} No funding was received for conducting this study.

\end{document}